\documentclass{amsart}
\newcommand{\Li}{\operatorname{Li}}
\newcommand{\Lip}{\operatorname{Li}^{(p)}}
\newcommand{\Lit}{\widetilde{\operatorname{Li}}}
\newcommand{\fli}{\operatorname{li}}
\newcommand{\C}{\mathbb{C}}
\newcommand{\Fbar}{\bar{\mathbb{F}}}
\newcommand{\Z}{\mathbb{Z}}
\newcommand{\Q}{\mathbb{Q}}

\newtheorem{theorem}{Theorem}[section]

\newtheorem{proposition}[theorem]{Proposition}
\newtheorem{lemma}[theorem]{Lemma}
\newtheorem{corollary}[theorem]{Corollary}
\theoremstyle{definition}

\newtheorem{remark}[theorem]{Remark}

\numberwithin{equation}{section}
\begin{document}
\title{Finite and $p$-adic polylogarithms}
\author{Amnon Besser}
%\date{8.5.2000}
\maketitle

\section{Introduction}
\label{sec:intro}
The finite logarithm was introduced by Kontsevich (under the name
``The $1\frac{1}{2}$ logarithm'') in~\cite{Kon99}. The finite logarithm is
the case $n=1$ of the $n$-th polylogarithm  $\fli_n\in \Z/p[z]$
defined by $\fli_n(z)= \sum_{k=1}^{p-1} z^k/k^n$. In loc.\ cit.
Kontsevich proved that the finite logarithm satisfies a 4-term
functional equation, known as the fundamental equation of information
theory. The same functional equation is satisfied by the so called
infinitesimal dilogarithm
$-(x\log|x|+(1-x)\log|1-x|)$. Cathelineau~\cite{Cat96}
defined general infinitesimal polylogarithms and found that they
satisfy interesting functional equations. It was the idea of
Elbaz-Vincent
and Gangl \cite{Elb-Gan99} that these functional equations should be satisfied
by finite polylogarithms (The name ``finite polylogarithm'' is due to
them). Inspired by their work, Kontsevich raised the idea that the
finite polylogarithm could be a reduction of an infinitesimal version
of the $p$-adic 
polylogarithm, as defined by Coleman~\cite{Col82}. If such a connection is
established, it makes sense to hope that functional equations of the
infinitesimal $p$-adic polylogarithm can be established in a similar
way to its complex counterpart and that these then imply by reduction
the functional equations of the finite polylogs. A conjectural formula
for the precise $p$-adic polylog whose ``derivative'', in the sense to
be explained below, reduces to the finite polylog
was formulated by Kontsevich and proved by him for small $n$.
The purpose of this short note is to prove such a connection between
$p$-adic polylogarithms and finite polylogarithms.

To state the main result we recall that Coleman defined $p$-adic
polylogarithms, $\Li_n:\C_p \to \C_p$. These functions are locally
analytic in the sense that they are given by a convergent power series
on each residue disc in $\C_p$. We define the differential operator
$D$ by $D=z(1-z) d/dz$. Let $\Fbar_p$ be the algebraic closure of the
finite field with $p$ elements and let $W=W(\Fbar_p)$ be the ring of
Witt vectors of $\Fbar_p$, so $W$ is the ring of integers of the
maximal unramified extension
of $\Q_p$. Let $\sigma:\Fbar_p \to \Fbar_p$ be the automorphism which
is the inverse of the $p$-power map. Let $X=\{z\in W:\quad
|z|=|z-1|=1\}$. Our main result is then:
\begin{theorem}\label{thetheorem}
  For every $n>1$ let
  \begin{equation*}
    F_n(z) = \sum_{k=0}^{n-1} a_k \log^k(z) \Li_{n-k}(z)\;,
  \end{equation*}
  with $a_0=-n$ and 
  \begin{equation*}
    a_k= \frac{(-1)^{k}}{(k-1)!} +\frac{(-1)^{k+1}n}{k!}\;,
  \end{equation*}
  for $k>0$. Then the following holds
  for every $p>n+1$: One has $DF_n(X) \subset p^{n-1} W$ and for
  every $z\in X$ one has
  $p^{1-n} DF_n(z) \equiv \fli_{n-1}(\sigma(z)) \pmod{p}$.
  Furthermore, the choice of the coefficients $a_k$ is the unique
  choice of coefficients in $\Q$ for which the theorem holds
  for all $p>n+1$.
\end{theorem}
We remark that for a given $p$ there will be many other choices of
coefficients, for example those which are sufficiently congruent to
the $a_k$.

The author would like to thank Kontsevich for explaining the
conjecture to him and to Gangl for various comments and corrections.
He would also like to thank Nekov\'a\v r who first suggested to him the
idea of a connection between finite and $p$-adic polylogarithms.

\section{The proof}
\label{sec:proof}

The connection between $p$-adic and finite polylogarithms is made by
the following
\begin{proposition}
  Let $\Lip_n(z):= \Li_n(z)- \Li_n(z^p)/p^n$. Then $\Lip_n(X) \subset
  W$ and the function $\Lip_n$ reduces modulo $p$ to $(1-z^p)^{-1}
  \fli_n(z)$.
\end{proposition}
\begin{proof}
According to \cite{Col82}, the function $\Lip_n(z)$ can be computed as
\begin{equation*}
  \Lip_n(z)= \int_{\Z_p^\times}x^{-n} d\mu_z(x)\;,
\end{equation*}
where $\mu_z$ is the measure on $\Z_p$ defined by
\begin{equation*}
  \mu_z(a+p^m \Z_p) = \frac{z^a}{1-z^{p^m}}\quad a=0,1,\ldots,p^m-1\;.
\end{equation*}
Since for $z\in X$ $\mu_z$ takes integral values, this shows the first
statement. Reducing modulo $p$ we may replace the 
function $x\mapsto x^{-k}$ by the function $x\mapsto a^{-k}$ if $x\equiv a
\pmod{p}$, which is congruent to it modulo $p$ on $\Z_p^\times$. This
implies that $\Lip_n(z)$ is congruent modulo $p$ to
\begin{equation*}
  \sum_{a=1}^{p-1} a^{-n} \mu_z(a+p\Z_p)=
  \sum_{a=1}^{p-1} a^{-n} \frac{z^a}{1-z^p}\;.
\end{equation*}
\end{proof}
\begin{corollary}\label{reduction}
  Let $\alpha\in X$ be a root of unity. Then we have
  $\Li_n(\alpha)\in p^n W$ and $p^{-n} \Li_n(\alpha) \equiv
  -\fli_n(\sigma(\alpha))/(1-\alpha) \pmod{p}$.
\end{corollary}
\begin{proof}
Since $\alpha\in X\subset W$ the order of $\alpha$ is prime to $p$ and
we have $\alpha^{p^k}=\alpha$ for some $k$. By using the
definition of $\Lip_n$ repeatedly we find
\begin{align*}
  \Li_n(\alpha)&=\Lip_n(\alpha)+p^{-n} \Li_n(\alpha^p)\\
               &=\Lip_n(\alpha)+p^{-n}
               \Lip_n(\alpha^p)+p^{-2n} \Li_n(\alpha^{p^2})\\
        \cdots &=\sum_{i=0}^{k-1} p^{-in}\Lip_n(\alpha^{p^i})
                  + p^{-kn} \Li_n(\alpha^{p^k})\;.
\end{align*}
Since $\alpha^{p^k}=\alpha$ we may move the last term to the
left hand side of the equation and obtain
\begin{align*}
  \Li_n(\alpha)&= \frac{1}{1-p^{-kn}}\sum_{i=0}^{k-1}
                  p^{-in}\Lip_n(\alpha^{p^i})\\
               &= \frac{p^n}{p^{kn}-1}\sum_{i=0}^{k-1}
  p^{(k-1-i)n}\Lip_n(\alpha^{p^i})\in p^n W \;,
\end{align*}
and dividing by $p^n$ and reducing modulo $p$ we obtain using the
proposition
\begin{equation*}
  - \Lip_n(\alpha^{p^{k-1}}) \equiv (1-(\alpha^{p^{k-1}})^p)^{-1}
  \fli_n(\alpha^{p^{k-1}}) \equiv -\fli_n(\sigma(\alpha))/(1-\alpha)
  \pmod{p}\;.
\end{equation*}
\end{proof}

\begin{proposition}\label{maincong}
  Let $\alpha$ be a root of unity in $X$. Set
  $\Lit_n(\alpha)= p^{-n} \Li_n(\alpha)$. Then for $w\in W$ we have
  \begin{equation*}
    p^{-n} \Li_n(\alpha(1+pw))\equiv \sum_{k=0}^n
    \Lit_{n-k}(\alpha)\frac{w^k}{k!}\pmod{p}\;.
  \end{equation*}
\end{proposition}
\begin{proof}
Let $g_n(w)=p^{-n} \Li_n(\alpha(1+pw))$. Then one finds
\begin{equation*}
  g_0(w)=\frac{\alpha(1+pw)}{1-\alpha(1+pw)}=
  \frac{\alpha}{1-\alpha}\frac{1+pw} {1-\frac{\alpha}{1-\alpha}pw}
  =\Lit_0(\alpha)+\sum_{k=1}^\infty b_k (pw)^k
\end{equation*}
with $b_n \in W$. Write
\begin{equation*}
  g_n(w)= \Lit_n(\alpha)+\sum_{k=1}^\infty d_k^n w^k\;.
\end{equation*}
It is easy to verify that
\begin{equation*}
  \frac{d}{dw} g_n(w)= g_{n-1}(w)\frac{1}{1+pw}=
  g_{n-1}(w)(1+\sum_{k=1}^\infty c_k (pw)^k)\; ,
\end{equation*}
with $c_k \in W$.
To find the coefficients $d_k^n$ modulo $p$ for $k<p$ one can simply reduce
the above equations modulo $p$ and one easily finds that
\begin{equation*}
  d_k^n \equiv
  \begin{cases}
    \frac{\Lit_{n-k}(\alpha)}{k!} &\text{ when } k\le n \\
        0 & \text{ when } n<k<p\;.
  \end{cases}
\end{equation*}
It thus remains to show that also for $k\ge p$ the coefficient $d_k^n$ is
divisible by $p$. For this it is easier to consider the function
$f_n(u)=\Li_n(\alpha+u)$ which satisfies $f_0(u)\in W[[u]]$, $d/du
f_{n+1}(u) = g(u)f_n(u)$ with $g(u)\in W[[u]]$ and $f_n(0)\in
W$.
\begin{lemma}
  In the situation above we have $v_p(a_k) \ge -v_p(k!)$, where $v_p$
  is the $p$-adic valuation and $a_k$ is the $k$-th coefficient in the
  power series expansion with respect to $u$ of any
  of the functions $f_n$.
\end{lemma}
\begin{proof}
Let $a_k(h)$ be the $k$-th coefficient of $h$ for any power series
$h$. We have 
\begin{equation*}
  v_p(a_k(f_n g)) \ge \min_{l\le k}[ v_p(a_l(f_n))]\;.
\end{equation*}
This implies that
\begin{equation*}
  v_p(a_k(f_{n+1})) \ge \min_{l< k} [v_p(a_l(f_n))] -v_p(k)\;.
\end{equation*}

The lemma is clearly true for $n=0$. Suppose it is true for $n$. Then
\begin{equation*}
  v_p(a_k(f_{n+1})) \ge \min_{l\le k} [-v_p(l!)] -v_p(k) \ge -v_p((k-1)!)
  -v_p(k)= -v_p(k!)\;.
\end{equation*}
\end{proof}
Since $g_n(w) = p^{-n} f_n(\alpha pw)$, to finish the proof we have to check
that for every $k\ge p$ we have $k-n-v_p(k!) >0$. The well known
estimate $v_p(k!)\le k/(p-1)$ and the assumption $n<p$ imply that it
is sufficient to require $k(1-1/(p-1))\ge p$, and this is satisfied
for $k\ge p+1$. The case $k=p$ is also OK according to this analysis
unless $n=p-1$. 
\end{proof}
\begin{proof}[Proof of Theorem~\ref{thetheorem}]
Using the fact that $D$ is a derivation and that
\begin{equation*}
  D \log^k(z) = z(1-z) \frac{d}{dz} \log^k(z)=z(1-z) k
  \log^{k-1}(z)\frac{1}{z} = (1-z) k \log^{k-1}(z)
\end{equation*}
and
\begin{equation*}
  D \Li_k(z) = z(1-z) \frac{d}{dz} \Li_k(z) = z(1-z) \Li_{k-1}(z)
  \frac{1}{z} = (1-z) \Li_{k-1}(z)
\end{equation*}
we see that
\begin{align*}
  DF_n(z)&=D \sum_{k=0}^{n-1} a_k \log^k(z) \Li_{n-k}(z)\\&=
  (1-z) \sum_{k=0}^{n-1} a_k(k \log^{k-1}(z) \Li_{n-k}(z) + \log^k(z)
  \Li_{n-k-1}(z)) \\&=
  (1-z) \sum_{k=0}^{n-1} \log^k(z) \Li_{n-k-1}(z) (a_k + (k+1) a_{k+1})\;.
\end{align*}
Here we understand that $a_n=0$. Every $z\in X$ can be written as
$\alpha (1+pw)$ with $\alpha$ a root of
unity in $X$ and  $w\in W$. We have 
\begin{equation*}
  p^{-1} \log(\alpha(1+pw)) = p^{-1} \log(1+pw) \equiv w \pmod{p}\;.
\end{equation*}
If we assume that $a_k\in Z[1/(n+1)!]$ we now
find from the last computation and from proposition~\ref{maincong},
\begin{align*}
  p^{1-n} &DF_n(\alpha(1+pw))\\ &\equiv 
  (1-\alpha) \sum_{k=0}^{n-1} w^k(a_k+ (k+1)
  a_{k+1})\sum_{m=0}^{n-k-1} \frac{1}{m!} \Lit_{n-k-1-m}(\alpha) w^m \\
  &\equiv
  (1-\alpha) \sum_{l=0}^{n-1} w^l \sum_{k=0}^l (a_k + (k+1) a_{k+1})
  \frac{1}{(l-k)!} \Lit_{n-l-1}(\alpha)\;.
\end{align*}
It follows that to make the reduction of $p^{1-n} DF_n$ independent of
$w$ for all $p>n+1$ it is necessary and sufficient that for
$l=1,\ldots,n-1$ we have
\begin{equation}\label{conds}
  \sum_{k=0}^l (a_k + (k+1) a_{k+1})
  \frac{1}{(l-k)!} =0\;.
\end{equation}
If this is satisfied then the reduction of $p^{1-n} DF_n(\alpha(1+pw))$ is 
$(a_0+a_1) (1-\alpha) \Lit_{n-1}(\alpha)\equiv -(a_0+a_1)
\fli_{n-1}(\sigma(\alpha))$ so we should also require $a_0+a_1=-1$.

Let $A(t) = \sum_{k=0}^{n-1} a_k t^k$. Then the relations
\eqref{conds} can be
written as $e^t (A(t) + dA(t)/dt) \equiv a \pmod{t^n}$ where $a$ is a
constant. This implies that $A(t)+dA(t)/dt \equiv ae^{-t} \pmod{t^n}$
and after solving the resulting differential equation that the $n-2$
first equations in~\eqref{conds} are equivalent to
$A(t) \equiv (at+b) e^{-t} \pmod{t^{n-1}}$ for some other constant
$b$. We have $-1=a_0+a_1=b+(a-b)=a$. Now, in $(b-t)e^{-t}$ the
coefficient of $t^n$ is
\begin{equation*}
   b \frac{(-1)^{n}}{n!} -\frac{(-1)^{n-1}}{(n-1)!} =
  \frac{(-1)^{n-1}}{(n-1)!}(1-b/n)\;.
\end{equation*}
It is easy to see that for the equation~\eqref{conds} with $l=n-1$ to be
satisfied, this coefficient must be $0$ and hence $b=n$. This gives
the choice of the coefficients in the theorem and shows that they are
the unique choice in $\Z[1/(n+1)!]$. Now if we have coefficients in
$\Q$ satisfying the theorem, then we may clear denominators not
dividing $(n+1)!$ and using only independence of $w$ we obtain that
these must be a rational multiple of our $a_k$. But since the
reduction of $p^{1-n} DF_n$ is non-trivial the multiplier must be $1$.
\end{proof}
\begin{remark}
We have the equation $F_n(z) =-n L_n(z) - L_{n-1}(z) \log(z)$, where
\begin{equation*}
L_n(z)= \sum_{m=0}^{n-1}\frac{(-1)^m}{m!}\Li_{n-m}(z) \log^m(z)
\end{equation*}
is the
function defined in \cite{Bes-deJ98}. By loc.\ cit.\ the 
function $F_n$ satisfies $F_n(z) + (-1)^n
F_n(1/z)=0$. Differentiating this relation one gets the relation
$z DF_n(1/z)+ (-1)^n DF_n(z)=0$. Reducing modulo $p$ we find
$z \fli_{n-1}(1/z) + (-1)^n \fli_{n-1}(z) = 0$. This relation is
easily verified directly.
\end{remark}

\section{Another proof of the main result}
\label{sec:another}

In this section we sketch another proof of the main result. This proof
has two interesting features: First of all, it proves directly the
formula for $F_n$ in terms of the functions $L_n(z)$ defined at the
end of the last section. This formula is of course simpler than
the original formula. The other feature is that the key ingredient in
the proof is a formula, discovered by Rob de Jeu and the author, which
seems to be of some further importance. This formula shows up in the
computation of syntomic regulators. These two features suggest that
the proof to be described below may in some way be more ``correct'' than
the first one, although it is if anything slightly more complicated.

The formula alluded to above is the content of the following
\newcommand{\dlog}{d\log}
\begin{proposition}
  We define a sequence of functions $f_k(z,S)$ inductively as follows:
  \begin{equation*}
    f_0(z,S) = \frac{S}{1-S},\quad f_{k+1}(z,S)= \int_z^S
    f_k(z,t)\dlog t\;.
  \end{equation*}
  Then, when $z,S\in W$ and $z\equiv S \pmod{p}$  the following
  formula holds:
  \begin{equation}\label{Delprop}
    -\sum _{k=0}^{n}(-1)^{k}k!\cdot \binom{n}{k}f_{k+1}(z,S)\log
    ^{n-k}(S)
    = (-1)^{n} n! \left( L_{n+1}(z) - L_{n+1}(S) \right)\;.
  \end{equation}
\end{proposition}
We would like to remark on the potential importance of this
formula. The construction of $p$-adic polylogarithms by Coleman is an
inductive procedure. At each step the degree $n$ polylogarithm $\Li_n$
is constructed as a locally analytic function satisfying the
differential equation $d \Li_n(z) = \Li_{n-1}(z) \dlog z$. This
determines $\Li_n$ up to a locally constant function and a Frobenius
condition replaces this by a globally constant function ambiguity. As
remarked by Kontsevich the distribution relation removes the ambiguity
completely. The formula above allows for a different approach: The
functions $f_n$ have no ambiguity in their definitions and are in fact
given by converging power series in $S$ and $z-S$. Once the $f_n$ are
given, the formula determines $L_n$ up to a locally constant function
and the distribution relation determines it completely. The functions
$\Li_n$ can be determined from the $L_n$.

For our purposes, the formula is also useful because it allows us to
relate the values of $L_n$ at two congruent points.
\begin{lemma}
  Suppose $z\in X$ and $w\in W$. Then
\begin{equation*}
  p^{-n}f_n(z,z(1+pw))\equiv \frac{z}{1-z}\frac{w^n}{n!} \pmod{p}\;.
\end{equation*}
\end{lemma}
The proof is a direct computation similar to the proof of
proposition~\ref{maincong}. Suppose now that $z$ is a
root of unity. Then $\log(z(1+pw))\equiv pw \pmod{p^2}$. Thus, using
\eqref{Delprop} we immediately obtain
\begin{equation*}
  p^{-n-1}(-1)^n n! (L_{n+1}(z(1+pw))-L_{n+1}(z))\equiv
  c_{n+1} \frac{z}{1-z}w^{n+1}\pmod{p}
\end{equation*}
with
\begin{equation*}
  c_{n+1}=\sum_{k=0}^{n} (-1)^k \frac{k!}{(k+1)!}\binom{n}{k}=
  \frac{1}{n+1}\sum_{k=0}^{n} (-1)^k \binom{n+1}{k+1}=\frac{1}{n+1}\;.
\end{equation*}
so
\begin{equation*}
    p^{-n} (L_n(z(1+pw))-L_n(z))\equiv -(-1)^n \frac{1}{n!}
    \frac{z}{1-z}w^n\pmod{p}\;.
\end{equation*}
We can define the derivation $D$ on functions of two variables as
$D=S(1-S)\partial/\partial S+ z(1-z)\partial/\partial z$. One easily obtains
the following
\begin{lemma}
  If $z,S\in X$ and $z\equiv S\pmod{p}$ then $p^k | Df_k(z,S)$.
\end{lemma}
Differentiating the key formula we find, with $S=z(1+pw)$,
\begin{multline*}
  p^{-n}(-1)^n n! (DL_{n+1}(S)-DL_{n+1}(z))\\=
  p^{-n} \sum_{k=0}^n k! \binom{n}{k}\left[Df_{k+1}(z,S)\log^{n-k}(S)+(n-k)
  f_{k+1}(z,S)\log^{n-k-1}(S) \right](1-S)\\ \equiv
  \sum_{k=0}^n k! \binom{n}{k}(n-k) f_{k+1}(z,S)\log^{n-k-1}(S) (1-z)
\equiv
  d_{n+1}zw^{n}\pmod{p}\; ,
\end{multline*}
where
\begin{equation*}
  d_{n+1}=\sum_{k=0}^{n} (-1)^k \frac{k!}{(k+1)!}\binom{n}{k}(n-k)=
  -\sum_{l=1}^{n} (-1)^l \binom{n}{l}=1 \;.
\end{equation*}
Let $z$ be a root of unity. Since $\log(z)=0$ we have
\begin{align*}
  DL_n(z)&=
  z(1-z)\Big[
             \Li_{n-1}(z)\frac{1}{z}\\
             &+\sum_{m=1}^{n-1} \frac{(-1)^m}{m!}
               \left(
                     \Li_{n-m-1}(z) \frac{1}{z} \log^m(z)
                      + m \Li_{n-m}(z) \log^{m-1}(z) \frac{1}{z}
              \right)
       \Big]\\
       &=(1-z) (\Li_{n-1}(z)-\Li_{n-1}(z))=0\;.
\end{align*}
Thus we find
\begin{equation*}
    p^{-n} DL_n(S)\equiv -(-1)^n \frac{n}{n!} zw^{n-1}\pmod{p}\;.
\end{equation*}
Suppose now that we let $F(z)=\sum_{m=1}^n e_m L_m(z)\log^{n-m}(z)$.
Substituting first of all a root of unity in $X$ we
find $DF(z)=e_{n-1} (1-z) L_{n-1}(z)$ so we should have $e_{n-1}=-1$.
Then
\begin{align*}
  p^{1-n}DF(S)&= \sum_{m=1}^n  e_m\Big(DL_m(S)\log^{n-m}(S)
  \\&\phantom{=} \phantom{\sum_{m=1}^n 
    e_m}+(n-m)L_m(S)\log^{n-m-1}(S)(1-S) \Big) p^{1-n}
  \\ &\equiv
  \sum_{m=1}^n e_m\Big(-(-1)^m \frac{m}{m!} zw^{m-1}\cdot
  w^{n-m}\\&\phantom{=} \phantom{\sum_{m=1}^n
    e_m} +
  (n-m)(p^{-m}L_m(z)-(-1)^m \frac{1}{m!}\frac{z}{1-z}w^m)w^{n-m-1}(1-z)\Big)
  \\ &=
  -zw^{n-1}\sum_{m=1}^n ((-1)^m e_m (\frac{m}{m!} + (n-m)\frac{1}{m!}))\\
  &\phantom{=} \phantom{-zw^{n-1}}+\sum_{m=1}^n e_m (n-m) p^{-m}L_m(z)
  \pmod{p} \;. 
\end{align*}
\newcommand{\notequiv}{not \equiv}
We wish to choose the coefficients $e_m$ in such a way that this expression
is independent of $w$. Since in the second sum we know by
corollary~\ref{reduction} that
$p^{-m}L_m(z)$ is not congruent to $0$ modulo $p$ (as a function of
$z$) we see that the only non zero
coefficients can be $e_n$ and $e_{n-1}$. By fixing $e_{n-1}=-1$ we get the
equation
\begin{equation*}
  (-1)^n e_n \frac{n}{n!}+(-1)^{n}
  \left(\frac{n-1}{(n-1)!}+\frac{1}{(n-1)!}\right)=0
\end{equation*}
from which we can recover
$e_{n}=-n$.

\end{document}